
\documentstyle{amsppt} 
\input amstex
\magnification=\magstep1
\TagsOnRight
\hsize=5.2in                                                  
\vsize=7.9in

\topmatter

\title Absolute torsion\endtitle
\dedicatory Dedicated to Mel Rothenberg on the occasion of his 65-th birthday\enddedicatory
\rightheadtext{Absolute torsion}
\leftheadtext{M. Farber and V. Turaev}
\author  Michael Farber* and Vladimir Turaev** \endauthor
\address
School of Mathematical Sciences,
Tel-Aviv University,
Ramat-Aviv 69978, Israel\endaddress
\address
 Institut de Recherche Math\'ematique Avanc\'ee,
Universit\'e Louis Pasteur - C.N.R.S., 7 rue Ren\'e Descartes, 67084 Strasbourg, France
\endaddress
\email farber\@math.tau.ac.il, turaev\@math.u-strasbg.fr
\endemail
\thanks{\footnotemark"*" Partially supported by a grant from the 
US - Israel Binational Science Foundation and by
the Herman Minkowski Center for Geometry\newline 
\footnotemark"**"
Partially supported by the EC TMR network "Algebraic Lie Representations",
EC-contract no ERB FMRX-CT97-0100
} 
\endthanks
\abstract{In this paper we use the results of our previous work \cite{FT}
in order to compute {\it the phase of the torsion} of an Euler structure $\xi$ in terms of 
the characteristic class $c(\xi)$. 
Also, we introduce here a new notion of {\it an absolute torsion}, which does not require 
a choice of any additional topological information (like an Euler structure). 
We prove that in the case 
of closed 3-manifolds obtained by 0-surgery on a knot in $S^3$ the 
absolute torsion is equivalent to the Conway polynomial. Hence the absolute torsion can be viewed as a
high-dimensional generalization of the Conway polynomial.}
\endabstract
\endtopmatter

\def\buildrel#1\over#2{\mathrel{\mathop{\null#2}\limits^{#1}}}
\def\buildrul#1\under#2{\mathrel{\mathop{\null#2}\limits_{#1}}}

\define\C{{\bold C}}
\define\R{{\bold R}}  
       
\define\Z{{\bold Z}}  
\define\Q{{\bold Q}}

\define\T{{\Cal T}}

\define\Hom{\operatorname{Hom}}

\redefine\det{\operatorname{det}}

\redefine\Im{\operatorname{Im}}   


\def\<{\langle}
\def\>{\rangle}
\def\kk{\bold k}

 \def\det  {{\text {det}}}

 \def\mod {{\text {mod}}}
 \def\Hom {{\text {Hom}}}
 \def\Im {{\text {Im}}}
 \def\dim   {{\text {dim}}}
 
 \def\Eul {{\text {Eul}}}
 
\define\Det{{\frak {Det}}}

\heading{\bf \S 1. Euler structures and Poincar\'e-Reidemeister metric}\endheading

In this section we give a brief review of the main results of \cite{FT}, which we will use in the sequel. 
The proofs of all theorems, appearing in this section, can be found in \cite{FT}.

\subheading{1.1. Determinant lines}
 We shall denote by $\kk$ a fixed ground field of characteristic zero.
The most important special cases are $\kk =\R$ and  $\kk=\C$.

If $V$ is a finite dimensional vector space over $\kk$, {\it the
determinant line of} $V$ is denoted by 
$\det \, V$ and is defined as the top
exterior power of $V$, i.e., $\Lambda^n V$, where $n = \dim \, V$. 
The dual line $\Hom_{\kk} (V,\kk)$ is denoted by $(\det \, V)^{-1}$.
For a finite dimensional graded vector space
$V = V_0\oplus V_1\oplus   \dots \oplus V_m $, its {\it determinant
line}  $\det \, V$ is defined as the tensor product 
$\det \, V = \det \, V_0 \otimes (\det \, V_1)^{-1}\otimes \det \, V_2
\otimes
\dots \otimes (\det \, V_m)^{(-1)^m}.$

Let $C$ be a finite
dimensional
chain complex 
over $\kk$.
In the theory of  torsion a crucial role is played by a canonical
isomorphism
$$\varphi_{C}: \det \,  C  \to \det \,  H_\ast(C), \tag1-1$$
where    both $C$
and $H_\ast(C)$ are considered   as graded vector spaces. 
The definition of the
mapping $\varphi_{C}$ is as follows.
 Choose for each $q=0,...,m$  
non-zero elements
$c_q\in \det \, C_q $ and $h_q\in \det \, H_q(C)$. Set 
$c=c_0\otimes c_1^{-1}\otimes c_2\otimes \dots \otimes c_m^{(-1)^m}\in
\det \, C $ and  $h=h_0\otimes h_1^{-1}\otimes h_2\otimes \dots
\otimes h_m^{(-1)^m} \in \det \,H_\ast(C)$,
 where   $-1$ in the exponent 
denotes the dual functional.
We define $\varphi_{C}$ by  
$\varphi_C (c) = (-1)^{N(C)}\, [c:h]\,h, $
where $N(C)$ is a residue modulo 2 defined  below and $ [c:h]$ is a
nonzero element of $\kk$, defined by 
$$ [c:h] =  \prod_{q=0}^m
 [d(b_{q+1})\hat h_q b_q/\hat c_q]^{(-1)^{ q+1 }}.\tag1-2$$
Here $b_q$ is a sequence of vectors of $C_q$ whose
 image $d(b_q)$ under the boundary homomorphism $d:C_q\to C_{q-1}$
is a basis of $\Im\, d$; the symbol $\hat h_q$ denotes a sequence of
cycles in $C_q$ such that the wedge product of their homology classes 
equals $h_q$;
the symbol $\hat c_q$ denotes a basis of
  $C_q$ whose wedge product  
equals $c_q$; the number
$ [d(b_{q+1})\hat h_q b_q/\hat  c_q]$ is the determinant of the matrix
 transforming $\hat  c_q$ into the 
basis $d(b_{q+1})\hat h_q b_q$ of $C_q$.
 The residue $N(C)$  is defined by
$$N(C)=\sum_{q=0}^m\alpha_q(C)\beta_q(C) \,(\mod\,2), \tag1-3$$
where
$$\alpha_q(C) = \sum_{j=0}^q \dim \, C_j\,(\mod\,2),
\quad  \beta_q(C) = \sum_{j=0}^q \dim \, H_j(C)\,(\mod\,2).\tag1-4$$

Formula (1-2) involves the  sign
 refinement of the standard formula, introduced in \cite{T1}.
We will  introduce now more sign involving factors
 in other natural maps arising in this setting. 

\subheading{1.2. The fusion isomorphism} For two  finite-dimensional graded vector spaces $V= V_0\oplus
V_1\oplus   \dots \oplus V_m$
 and $W= W_0\oplus W_1\oplus   \dots \oplus W_m$, we
define a canonical isomorphism  
$$\mu_{V,W}:  \det \, V\otimes \det \, W \to    \det  (V\oplus
W), \tag1-5$$
  by  
$$\mu_{V,W}= (-1)^{M(V,W)}  \bigotimes_q \mu_q^{(-1)^q}, \tag1-6$$
where 
$\mu_q: \det \, V_q\otimes \det \, W_q \to    \det \,
(V_q\oplus W_q)$ is the isomorphism defined by
  $$(v_1\wedge v_2\wedge \dots \wedge v_k )\otimes
(w_1\wedge w_2\wedge \dots \wedge w_l) \mapsto v_1\wedge v_2\wedge
\dots \wedge v_l \wedge w_1\wedge w_2\wedge \dots \wedge w_k,$$
with $k=\dim \, V_q, l=\dim \, W_q  $, the isomorphism  
$\mu_q^{-1}$
is defined as the transpose of the inverse  of $\mu_q$, and
$${M(V,W)}= \sum_{q=1}^m \alpha_{q-1}(V)\,\alpha_q(W) \,\, \in 
\Z/2\Z, \tag1-7$$ 
with
$\alpha_{q-1}(V)$ and  $\alpha_q(W)$ defined as in (1-4).

We will call (1-5) {\it the fusion homomorphism.}

\subheading{1.3. Duality operator} 
Let 
$V = V_0\oplus V_1\oplus   \dots \oplus V_m $
be a finite dimensional graded vector space over $\kk$ with odd   $m$.
We define the  {\it  dual graded vector space} over $\kk$
by $V'=  V'_0\oplus V'_1\oplus   \dots \oplus V'_m $
where
$V'_q=(V_{m-q})^*=\Hom_{\kk} (V_{m-q},\kk)$.
We define a {\it duality operator} 
$$D=D_V:\det\,V \to \det\,V'$$ 
as follows. Let $v_q\in \det \, V_q$ be a volume element determined 
by a basis of $V_q$ and let $v'_{m-q}\in \det \, V'_{m-q}$ be the volume
element determined by the dual basis of $V'_{m-q}$, for $q=0,1,...,m$. 
Then   
$$D(v_0\otimes v_1^{-1}\otimes v_2\otimes \dots \otimes v_m^{-1})
=(-1)^{s(V)} v'_0\otimes (v'_1)^{-1}\otimes v'_2\otimes \dots \otimes
(v'_m)^{-1},$$ where
the residue $s(V)\in \Z/2\Z$ is given by
$$s(V)=\sum_{q=1}^m\alpha_{q-1} (V)\, \alpha_q (V)
+\sum_{q=0}^{(m-1)/2}\alpha_{2q}(V) .\tag1-8$$ 

\subheading{1.4. Euler structures and their characteristic classes}
We recall the notion of combinatorial Euler structure
on a CW-space, following \cite {T2}.

Let $X$ be a finite
connected CW-space with $\chi(X)=0$.
 An {\it Euler chain in $X$} is a 
singular 1-chain $\xi$ in $X$  such that
$d\xi = \sum_a (-1)^{|a|}p_a$
where $a$ runs over all cells of $X$ and $p_a $
is a point in   $a$; the symbol $|a|$ denotes the
dimension of $a$.
 
An {\it Euler structure on
$X$} is an equivalence class of Euler chains with respect to a natural
equivalence relation. The set of equivalence classes (i.e., the set of Euler structures on
$X$) is denoted by $\Eul(X )$. We shall usually
 denote an Euler structure
and a representing  it Euler chain by the same letter.

The group $H_1(X)$ acts on the set $\Eul(X )$ and this action is free and transitive.
We shall use multiplicative notation both for this action and  for the
group operation in $H_1(X)$.

Assume that $X$ is a closed connected PL-manifold with $\chi(X)=0$.
For each  Euler structure $\xi$ on 
$X$ we   define its {\it  characteristic class}  $c(\xi)  \in
H_1(X)$ as follows. Choose a
PL-triangulation $\rho$ of $X$.  Let $W$ be the   1-chain in $X$
defined by $W= \sum_{a_0<a_1 \in \rho} (-1)^{|a_0|+|a_1|}
 \langle {\underline {a}}_0,  {\underline {a}}_1 \rangle,$
where $ a_1$ runs over all simplices of $\rho$, $a_0$ runs over all
proper faces of $a_1$, and 
$\langle {\underline {a}}_0,  {\underline {a}}_1 \rangle$ is a path
in $  a_1$
going from the barycenter ${\underline {a}}_0$
of $a_0$ to the barycenter ${\underline {a}}_1$
of $a_1$. It is easy to check (see \cite{HT})  that 
$\partial W=(1-(-1)^m) \sum_{a\in \rho } (-1)^{|a | }
  {\underline {a}} $
where $m=\dim\, X$. Now, any
Euler structure  
 on $ X $ can be presented by   an Euler
chain $\xi$  in
$ (X,\rho)$  such that $ \partial \xi= \sum_{a } (-1)^{|a | }
  {\underline {a}} $. It is clear that $(1-(-1)^m)\, \xi- W$ is a 1-cycle.
Denote its class in $H_1(X)$ by $c(\xi)$. 
In this way, we obtain  a
mapping  $c:\Eul (X ) \to H_1(X)$.

If $m$ is odd, then  (in multiplicative notation)
$c(h\xi) = h^2c(\xi) $
for any $\xi\in \Eul (X), h\in H_1(X)$. 
For odd $m$, the mod 2 reduction of $c(\xi)$ is independent of 
$\xi$
and equals to the dual of the Stiefel-Whitney class 
$w_{m-1}(X)\in H^{m-1}(X,\Z/2\Z)$. This follows from the
fact that   $W\,(\mod\, 2)$ represents the dual of
$w_{m-1}(X)$, see \cite {HT}.

\subheading{1.5. Torsion of Euler structures}
Let $F$ be a flat $\kk$-vector bundle   over a finite
connected CW-space $X$.  For each Euler structure   $\xi\in \Eul(X)$  on $X$ we define a torsion 
$\tau(X,\xi;F)$ which is an element of the determinant line
$\det \, H_\ast(X;F)$ defined up to multiplication by
$(-1)^{\dim\, F}$. We denote by $C_\ast(X;F)$ the cellular chain
complex computing the homology of $X$ with values in $F$. Recall that
$$C_q(X; F)  = \bigoplus_{\dim \, a=q}   \, \Gamma(a, F),\tag1-9$$ 
where $\Gamma(a,F)$ denotes the space of flat sections of $F$ over $a$.
Set $$\tau(X,\xi;F)=\varphi_{C}
(c_0\otimes c_1^{-1}\otimes c_2\otimes \dots \otimes c_m^{(-1)^m})
\in \det \,  H_\ast(X;F)\tag1-10 $$
where $m=\dim\, X$ and   $c_q\in  
\det \,C_q(X;F)\, (q=0,1,...,m)$ are non-zero elements defined as follows.
Fix a point   $x\in X$
and a basis $e_x$ in the fiber $F_x$.
Let
$\beta_a:[0,1]\to X$ be a path connecting $x=\beta_a(0)$ to a point   
$ \beta_a(1) \in   a$. The assumption
$\chi(X)=0$ implies   that the 1-chain
$\sum_{a } (-1)^{|a|}   \beta_a$ (where $a$ runs over all cells of $X$)
is an Euler chain   with boundary $\sum_{a } (-1)^{|a|}  
\beta_a(1)$. We choose the paths $\{\beta_a\}_a$ so that this chain
represents   $\xi$. 
We apply the parallel
transport  to $e_x$ along $\beta_a$ to obtain a
basis in the fiber   $F_{\beta_a(1)}$ and we
extend it   to a basis of flat sections over $a$. The concatenation of
these bases over all $q$-dimensional cells gives a basis  in
$C_q(X;F)$ via (1-9). The wedge
product of the elements of this basis yields  
  $c_q\in   \det \,C_q(X;F)$. 

\proclaim{1.6. Lemma} If $\dim \, F$ is even then the torsion 
$\tau(X,\xi;F)\in \det \, H_\ast(X;F)$ is well defined, has no indeterminacy
and is combinatorially invariant.
If $\dim \, F$ odd, the torsion $\tau(X,\xi;F)$ has a sign indeterminacy. 
\endproclaim

In order to fix the sign of the torsion in the case $\dim \, F$ odd, 
one may use technique of {\it homological orientations}, i.e., the orientations
of the determinant line $\det \, H_\ast(X,\R)$ of the real cohomology, which was introduced in
 \cite{T1, T2}. Cf. also 3.3.

\subheading{1.7. The Poincar\'e-Reidemeister scalar product} 
Let $X$ be a  {\it closed connected oriented
piecewise linear manifold of odd dimension m.} Let $F$ be a
flat $\kk$-vector bundle over  $X$.
The standard  homological intersection pairing 
$$H_q(X;F^*)\otimes H_{m-q}(X;F)\to \kk\tag1-11$$ 
allows us to identify  the dual of $H_{m-q}(X;F)$ with $H_q(X;F^*)$.
Applying
  the construction of Section 1.3 to 
the graded vector space $\oplus_{q=0}^m H_q(X;F)$ we obtain
  a canonical isomorphism $$D: \det \, H_\ast(X;F)\to \det \,
H_\ast(X;F^\ast).\tag1-12$$  
It is easy to check that $D$ does not depend on the choice of
the
orientation of $X$. 

{\it The Poincar\'e-Reidemeister scalar product} is defined as the bilinear pairing
$$\langle \ ,\ \rangle_{PR}: \det \, H_\ast(X;F)\times \det \, H_\ast(X,
F) \to \kk,\tag1-13$$ given by
$$\langle a ,b \rangle_{PR} =
  \mu (a\otimes D(b)) /  \tau(X;F\oplus F^\ast)  \in \kk,  $$ 
where 
$a,b\in 
\det \, H_\ast(X;F)$ and   $D$ is the  
isomorphism (1-12). 
Here $\mu $ denotes
  the canonical fusion isomorphism
$$ 
   \det \, H_\ast(X; F) \otimes
\det \, H_\ast(X;  F_\ast)\to \det \,
 (H_\ast(X; F)\oplus H_\ast(X; F^\ast)) 
= \det \, H_\ast(X; F\oplus F^\ast)
$$ 
  defined in Section 1.2.

The Poincar\'e-Reidemeister scalar product determines {\it the
Poincar\'e-Reidemeister metric}
(or norm)
on the determinant line $\det \, H_\ast(X;F)$, which was introduced in
\cite{Fa}.
It is given by
$$a\mapsto \sqrt{|\<a,a\>_{PR}|},\quad a\in \det \, H_\ast(X;F)$$
(the positive square root of the absolute value of $\<a,a\>_{PR}$).
The PR-scalar product contains an additional phase or sign information.

The following Theorem computes the PR-scalar product in terms of the Euler structures.

\proclaim{1.8. Theorem } 
Let     $F$ be a 
flat $\kk$-vector bundle over a closed connected orientable 
PL-manifold
$X$ of odd dimension $m$. If $\dim\, F$ is odd, then we additionally assume
that $X$ is provided with a homology orientation.
 Then for
  any   $\xi\in
\Eul(X)$,
$$\langle \tau (X,\xi; F) ,\tau (X, \xi; F) \rangle_{PR}
= (-1)^z \det_F(c(\xi)), \tag1-14$$  
where $\langle \, ,\, \rangle_{PR}$ is the 
Poincar\'e-Reidemeister
scalar product, $\det_F(c(\xi))\in \kk^\ast$ is the determinant of the monodromy of $F$ along the characteristic
class $c(\xi)\in H_1(X)$, and the residue
$z\in \Z/2\Z$ is
$$
z=\cases
0,~ {\text  {if}}\,\,\,\dim\,F\,\,\, is \,\,\, even \,\,\, or  \,\,\, 
m\equiv 3\, (\mod\, 4), \\
s\chi(X)\, (\mod\, 2),~ {\text  {if}} \,\,\,
\dim\,F\,\,\, is \,\,\, odd \,\,\, and  \,\,\,  m\equiv 1 \,(\mod\,
4), \endcases \tag 1-15$$ 
where
$s\chi(X) = \sum_{i=0}^{(m-1)/2} \dim H^{2i}(X;\R)$
denotes the {\it semi-characteristic} of $X$.
 \endproclaim

The following   theorem describes the sign of the PR-pairing.

 \proclaim{1.9. Theorem} Let  $F$ be a 
flat $\R$-vector bundle over a closed connected orientable
PL-manifold
$X$ of odd dimension $m$.
The Poincar\'e-Reidemeister scalar product on $\det \, H_\ast(X;F)$ is 
positive  definite for $m\equiv 3\, (\mod\, 4)$. If $m\equiv 1\, (\mod\, 4)$ then the sign of the
Poincar\'e-Reidemeister scalar product equals
$$(-1)^{\<w_1(F)\cup w_{m-1}(X),[X]\>+s\chi(X)\cdot \dim F}.\tag1-16$$
\endproclaim

 Note two interesting special cases of Theorem 1.9 assuming $m\equiv 1\, (\mod\, 4)$.
If $w_{m-1}(X)=0$ and $\dim\, F$ is even then the Poincar\'e-Reidemeister scalar product is
positive definite.
The same conclusion holds if 
$F$ is orientable and $\dim\,F$ is even.

\proclaim{1.10. Theorem (Analytic torsion and Euler structures)} Let 
$X$ be a closed connected orientable smooth manifold of odd dimension  
and let $F$ be a flat $\R$-vector bundle over $X$.
If $\dim\, F$ is odd, then we additionally assume
that $X$ is provided with a homology orientation.
For any Euler structure $\xi\in Eul(X)$,
 the Ray-Singer norm of   its cohomological torsion 
$\tau^\bullet (X, \xi; F)\in \det \, H^\ast(X; F)$ (defined similarly to 1.5)
 is equal to
the  positive  square root of the absolute value of the monodromy of $F$
along the characteristic 
class $c(\xi)\in H_1(X)$:
\endproclaim
$$||\tau^\bullet(X, \xi;F)||^{RS} = |\det_Fc(\xi)|^{1/2}.\tag1-17$$

In the special case, where the flat bundle $F$ is acyclic, i.e., 
 $H^\ast(X;F)=0$, the torsion 
$\tau^\bullet(X, \xi; F)$ is a real number 
 and Theorem 1.10
yields
$$\prod_{q=0}^{\dim X} (\Det \, \Delta'_q)^{(-1)^{q+1}\, q}= 
\frac{(\tau^\bullet(X, \xi; F))^2}{|\det_Fc(\xi)|}.\tag1-18$$

Theorem 1.10 generalizes both the classical Cheeger-M\"uller theorem \cite{C},
\cite{Mu1}  (dealing with orthogonal flat real bundles $F$)  and also
 the (more general) theorem of M\"uller \cite{Mu2}
(dealing with the unimodular flat real bundles $F$). 
Note that if $F$ is unimodular then   
$|\det_Fc(\xi)|=1$ and the torsion $\tau^\bullet(X, \xi;F)$ does not
 depend on the choice
of $\xi$.

 \heading{\bf \S 2. Phase of the torsion}
\endheading

The purpose of this section is to give a formula expressing the phase of the torsion of Euler
structures $\xi\in \Eul(X)$ 
(understood as an element of a determinant line) in terms of the characteristic class $c(\xi)\in H_1(X)$.

Throughout this section $\kk=\C$.

\subheading{2.1. Involution on the determinant line} Let $X$ be a closed orientable piecewise linear
manifold of odd dimension $m$. Let $F\to X$ be a flat complex vector bundle admitting a flat Hermitian metric. 

We introduce a canonical involution on the complex line $\det \, H_\ast(X;F)$.
Consider the duality operator (1-12).
Recall that $F^\ast$ denotes the dual flat vector bundle bundle. 
Any flat Hermitian metric on $F$ determines an anti-linear
isomorphism of flat bundles 
$F^\ast \to F,$ 
which induces an anti-linear isomorphism
$$\psi: \det \, H_\ast(X;F^\ast)\to \det \, H_\ast(X;F).\tag2-1$$
\proclaim{Definition} The canonical involution on the determinant line $\det \, H_\ast(X;F)$
is the following anti-linear isomorphism
$$\overset-\to{}: \det \, H_\ast(X;F)\to \det \, H_\ast(X;F),\tag2-2$$
where for $\tau \in \det \, H_\ast(X;F)$ we set
$$
\overset-\to\tau =
(-1)^{s\chi(X)\cdot\dim F\cdot(m+1)/2}\psi(D(\tau)),\tag2-3
$$
where $m=\dim X$.
\endproclaim
Here $s\chi(X)$ denotes the semicharacteristic of $X$, i.e. 
$s\chi(X) = \sum_{i=0}^{(m-1)/2} b_{2i}(X).$

\proclaim{2.2. Lemma} (A) The anti-linear isomorphism (2-2) is an involution.\newline
(B) It is independent of the choice of a Hermitian metric on $F$.\newline
(C) If the flat bundle $F$ is acyclic then the determinant line $\det \, H_\ast (X;F)$ is canonically isomorphic to $\C$ and under this isomorphism the involution (2-2) coincides with the complex
conjugation.\endproclaim

A proof is given below.

An element $\tau \in \det \, H_\ast(X;F)$ will be called {\it real} if $\overset -\to \tau =\tau$. The real
elements of $\det \, H_\ast(X;F)$ form a real line.

For a nonzero $\tau \in \det \, H_\ast(X;F)$, its {\it phase} is defined as angle $\phi\in \R$
so that $\tau$ can be represented in the form $\tau = \tau_0 e^{i\phi}$, where $\tau_0$ is real. It is clear
that the phase $\phi$ is determined up to adding integral multiples of $\pi$. 
We will denote it by $\bold {Ph}(\tau)$.

The following theorem computes the phase of the torsion in terms of the characteristic
class of the Euler structure.

\proclaim{2.3. Theorem} 
Let $X$ be a closed orientable piecewise linear
manifold of odd dimension $m$ and let $F\to X$ be a flat complex vector bundle admitting a flat Hermitian metric. Then
the phase of the torsion $\tau(\xi, F)$ is given by the following formula:
$$
{\bold {Ph}}(\tau(X,\xi; F))= 
\frac{1}{2}\arg\det_F (c(\xi))\, \, \, mod \, \pi\Z.\tag2-4
$$
\endproclaim
Note that in the case when $\dim \, F$ is odd, the torsion $\tau(X,\xi;F)$ is defined only up to a sign,
but its phase is still well defined. 

The number $d = \det_F (c(\xi))$ lies on the unit circle and $\arg(d) =\arg\det_F (c(\xi))\in \R/(2\pi\Z)$ 
is defined by $d=\exp(i\arg(d))$.

\demo{Proof of Lemma 2.2 and Theorem 2.3} Let $\tau$ denote $\tau(X,\xi;F)$ (in the case when
$\dim F$ is even) or 
$\tau(X,\eta, \xi;F)$ (in the case when
$\dim F$ is odd); here $\xi\in \Eul(X)$ is an Euler structure and 
$\eta$ is a homological orientation (i.e. an orientation of the line $\det\, H_\ast(X;\R)$)
which we need in order to fix the
sign of the torsion in the case when $\dim F$ is odd, cf. 6.3 of \cite{FT}.

We first show that
$$\overline \tau = \det_F(c(\xi))^{-1}\cdot \tau.\tag2-5$$
We will use theorem 7.2 of \cite{FT}, which states 
$$D(\tau) = (-1)^{\dim F\cdot s\chi(X)\cdot (m+1)/2}\cdot \tau(X,\eta,\xi^\ast;F^\ast),\tag2-6$$
where $\xi^\ast\in \Eul(X)$ is the dual Euler structure. Since $\xi=c(\xi)\xi^\ast$ (cf. formula
(5-4) in \cite{FT}) we obtain
$$
\aligned
D(\tau) = &(-1)^{\dim F\cdot s\chi(X)\cdot (m+1)/2}\cdot \det_{F^\ast}(c(\xi))^{-1}\cdot \tau(X,\eta,\xi;F^\ast)=\\
&(-1)^{\dim F\cdot s\chi(X)\cdot (m+1)/2}\cdot \det_{F}(c(\xi))\cdot \tau(X,\eta,\xi;F^\ast),
\endaligned
\tag2-7$$
where $m=\dim X$.
Applying to both sides of (2-7) the anti-linear isomorphism (2-1) and using our definition (2-3)
we obtain (2-5).

To prove statement (A) of Lemma 2.2, we note that according to (2-5) we have
$$\overset-\to \tau = e^{i\psi}\tau,\quad \tau\in \det\, H_\ast (X;F),\quad \tau\ne 0\tag2-8$$
for some angle $\psi\in \R$. Here $\tau=\tau(X,\eta, \xi;F)$ and $\det_F(c(\xi)) = e^{i\psi}$. 
Therefore, we obtain
$\overline{\overline\tau} = e^{-i\psi}{\overset-\to\tau} = \tau.$
This proves that the anti-linear isomorphism (2-2) is involutive on the torsion $\tau$ and therefore it is involutive on any other element. 

(B) obviously follows from (2-5) since the torsion $\tau$ does not depend on the Hermitian metric
on $F$.

To prove (C) we observe that in the acyclic case the duality isomorphism $D$ (cf. (1.3)) can be identified (after the canonical identification of the determinant lines $\det\, H_\ast(X;F)$ and 
$\det\, H_\ast(X;F^\ast)$
with $\C$) with the multiplication by
$$
(-1)^{s\chi(X)\cdot\dim F\cdot (m+1)/2}.
$$
The isomorphism $\psi: \det \, H_\ast(X;F^\ast)\to \det \, H_\ast(X;F)$ induced by (2-3) after these identifications coincides with the usual
complex conjugation (as one sees from the definition of torsion). This implies our statement.

Now we will prove Theorem 2.2. 
If $\tau= e^{i\phi}\tau_0$, where $\tau_0$ is real (with respect 
to the canonical involution (2-2)), then $\overset-\to\tau= e^{-i\phi}\tau_0$ and from (2-5) we obtain
$\overset-\to\tau/\tau = e^{-2i\phi}=e^{-i\arg\det_F (c(\xi))}$.
Therefore,
$${\bold {Ph}}(\tau(X,\xi; F))= \phi = \frac{1}{2}\arg\det_F (c(\xi))\, \, \mod \, \pi \Z.$$
\qed
\enddemo

\heading{\bf \S 3. The Absolute Torsion}\endheading

In this section we introduce a new concept of torsion which we call {\it absolute torsion}. 
It has some important advantages with respect to other similar notions of torsion: on one hand it 
is well defined and has no indeterminacy in most important cases including non-unimodular flat 
bundles. On the other hand it requires no additional topological information, such as Euler structures.
We will show in the next section that the absolute torsion can be viewed as a natural high dimensional 
generalization of the Conway polynomial.

\subheading{3.1. Basic assumption} In this section we will always deal with a closed {\it oriented}
PL manifolds $X$ of odd dimension $m$ and a flat complex vector bundle $F\to X$, 
satisfying the following condition:

{\it (i) The Stiefel-Whitney class $w_{m-1}(X)\in
H^{m-1}(X,\Z_2)$ vanishes;

(ii) The first Stiefel-Whitney class $w_1(F)$, viewed as a homomorphism $H_1(X;\Z)\to \Z_2$,
vanishes on the 2-torsion subgroup of  $H_1(X;\Z)$.}

Note that condition (i) is automatically satisfied in the case $m\equiv 3\, (\mod \, 4)$,
as proven by W. Massey \cite{Ma}. The condition (ii) holds for any orientable bundle $F$.
Also, (ii) holds for any $F$ assuming that $H_1(X)$ has no 2-torsion.

\subheading{3.2. Canonical Euler structures} An Euler structure $\xi\in \Eul(X)$ will be called 
{\it canonical} if it has trivial characteristic class $c(\xi)\, =\, 0 \in H_1(X)$. Since the mod 2 reduction
of $c(\xi)$ coincides with $w_{m-1}(X)$, the assumption 3.1.(i) is necessary for the 
existence of canonical Euler structures. It is also sufficient: if $w_{m-1}(X)=0$ 
then $c(\xi)$ admits a square
root $c(\xi)^{1/2}\in H_1(X)$ and $\xi' = c(\xi)^{-1/2}\cdot \xi$ is a canonical Euler structure.

Since $c(h\xi) = h^2c(\xi)$, where $h\in H_1(X)$, we see that {\it 
the canonical Euler structure is unique if and only if the group $H_1(X)$ has no 2-torsion. }

In general, the number of canonical Euler structures on $X$ equals the order 
of the 2-torsion subgroup of $H_1(X)$. 

\subheading{3.3. Canonical homological orientation} Recall that {\it a homological orientation} 
of $X$
is a choice of an orientation of the line $\det\, H_\ast(X;\R)$. 
It is observed in \cite{T1} that if $X$ is an oriented closed odd-dimensional manifold then 
there exists {\it a canonical homological orientation} of $X$, 
which sometimes depends on the choice of orientation of $X$.
It is described as follows.

Fix an orientation on $X$. We assume that the dimension $m=\dim \, X$ is odd $m=2r+1$. 
For any $i\le r$ fix an arbitrary basis $h^i_1, \dots, h^i_{b_i}\in H_i(X;\R)$ and let 
$h^{m-i}_1, \dots, h^{m-i}_{b_i}\in H_{m-i}(X;\R)$ be the dual basis, i.e. 
$$\<h_k^i \times h_l^{m-i},[X]\>=\delta_{k,l}$$
where $[X]\in H_m(X\times X, X\times X-\Delta;\R)$ is the class corresponding 
to the given orientation. Here $\Delta\subset X\times X$ denotes the diagonal. This gives
volume forms $h^i=h_1^i\wedge \dots \wedge h_{b_i}^i\in \det \, H_i(X;\R)$ and 
$h^{m-i}=h_1^{m-i}\wedge \dots \wedge h_{b_i}^{m-i}\in \det \, H_{m-i}(X;\R)$ and hence we obtain
a canonical nonzero element 
$$h=h^0\otimes (h^1)^{-1}\otimes h^2\otimes \dots \otimes (h^m)^{-1}\in \det \, H_\ast(X;\R).$$

Suppose that we reverse the orientation of $X$. Then the volume elements $h^i$ with $i\le r$ will
be the same and each $h^i$ with $i>r$ will be replaced by $(-1)^{b_i}h_i$. Hence we obtain that
reversing the orientation of $X$ changes $h$ as follows: $h\mapsto (-1)^{s\chi(X)}\cdot h$.
We arrive at the following:

\proclaim{3.4. Proposition} Any closed oriented 
odd-dimensional manifold $X$ has a canonical homological orientation. 
The canonical homological orientations of $X$ does not depend on the orientation of $X$ if
and only if $s\chi(X)$ is even.
\endproclaim

\subheading{3.5. Definition of absolute torsion}
Our purpose in this section is to define a combinatorial invariant
$$\T(F) \in \det \, H_\ast(X;F)\tag3-1$$
for {\it arbitrary} flat $\C$-bundle $F$ over a closed oriented odd-dimensional manifold $X$ and a
complex flat vector bundle $F\to X$, 
satisfying the 
conditions 3.1. 

We emphasize that {\it we do not require $F$ to be unimodular}. 

We construct $\T(F)$ as follows. Choose a canonical Euler structure $\xi\in
\Eul(X)$ and consider
the torsion 
$$\T(F) = \tau(X,\eta,\xi;F)  \in \det \, H_\ast(X;F),\tag3-2$$
where $\eta$ is the canonical homological orientation, cf. 3.3, 3.4.
The result {\it will not depend} on the choice of the canonical Euler
structure $\xi$ because of our assumption (ii) in 3.1. Indeed, replacing the canonical Euler structure
$\xi$ by another one, $h\xi$, where $h$ belongs to the 2-torsion subgroup of $H_1(X)$, gives the 
following torsion $\tau(X,\eta, h\xi;F)=\det_F(h)\cdot \tau(X,\eta,\xi;F)=(-1)^{\<w_1(F),h\>}\cdot \tau(X,\eta,\xi;F)$,
and our statement now follows from condition (ii) in 3.1.

We call $\T(F)$ {\it the absolute torsion.} 

An equivalent way to construct the absolute torsion consists in the following. Pick an arbitrary
Euler structure $\xi\in \Eul(X)$. 
Note that by our condition (i) in 3.1, {\it the characteristic class $c(\xi)\in H_1(X)$ is a
square.} In fact, (in the additive notations) the mod 2 reduction of $c(\xi)$ vanishes,
since it is Poincar\'e dual of the Stiefel-Whitney class
$w_{m-1}$, which we assume to be zero.
We know that there
exists a square root
$c(\xi)^{1/2}\in H_1(X).$
The indeterminacy in computing $c(\xi)^{1/2}$ can be described as follows
$h\mapsto h\cdot c(\xi)^{1/2}$, where $h$ belongs to 2-torsion of $H_1(X)$. By condition (ii)
in (3.1)
the monodromy of $F$ along any loop representing the class $c(\xi)^{-1/2}$,
is well defined and we may set
$$\T(F) = \det_F(c(\xi)^{-1/2})\cdot \tau(X,\xi;F).\tag3-3$$

\proclaim{3.6. Conclusion} Assume that $X$ and $F$ satisfy conditions 3.1. The absolute torsion $\T(F) \in \det \, H_\ast(X;F)$ is well defined. It is independent of the orientation of $X$ under any
of the following conditions:
\roster
\item  If $F$ is an even dimensional flat bundle; 
\item If $F$ is an odd dimensional flat bundle and the semi-characteristic $s\chi(X)$ is even.
\endroster
If $F$ is an odd dimensional flat bundle and the semi-characteristic $s\chi(X)$ is of $X$ odd, the absolute torsion changes the sign, when the orientation of $X$ is reversed.
\endproclaim

We now establish some properties of the absolute torsion.

\proclaim{3.7. Theorem (Duality)} Let $F$ be a flat complex vector bundle over $X$ and
let $F^\ast$ denote the dual flat vector bundle. Let
$$D: \det \, H_\ast(X;F)\to \det \, H_\ast(X;F^\ast).\tag3-4$$
be the duality operator (1-12). 
Then, 
$$D(\T(F))\, =\, (-1)^{s\chi(X)\cdot\dim F\cdot(m+1)/2}\T(F^\ast),\tag3-5$$
where $m=\dim X$.
\endproclaim
\demo{Proof} As
we showed in section \S 2 (cf. (2-7)), 
$$
\aligned
&D(\tau(X, \eta,\xi; F)) =\\
&(-1)^{s\chi(X)\cdot\dim F\cdot(m+1)/2}\cdot
\det_F(c(\xi))\cdot \tau(X,\eta,\xi; F^\ast)\in \det \, H_\ast(X; F^\ast).
\endaligned
\tag3-6$$
Dividing both sides by $\det_F(c(\xi)^{1/2})$ and observing that
$$\det_{F^\ast}(c(\xi)^{-1/2})= \det_F(c(\xi)^{1/2})$$
we obtain
$$D(\frac{\tau(X,\eta,\xi; F)}{\det_F(c(\xi)^{1/2})}) =
(-1)^{s\chi(X)\cdot\dim F\cdot(m+1)/2}\cdot
\frac{\tau(X,\eta,\xi; F^\ast)}{\det_{F^\ast}(c(\xi)^{1/2})}$$
which, combined with the definition, proves out statement. \qed
\enddemo

As the corollary we obtain that the absolute torsion is always real:

\proclaim{3.8. Theorem} Assume that $X$ and $F$ satisfy conditions 3.1.
Suppose that $F$ admits a flat Hermitian metric. Recall the canonical involution
on $\det \, H_\ast(X; F)$, cf. \S 2. Then 
$\T(F)\in \det \, H_\ast(X; F)$ is real:
$$\overline{\T(F)} = \T(F).\tag3-7$$
\endproclaim
\demo{Proof} Apply isomorphism (2-2), defined by a flat Hermitian metric on $X$, to both
sides of (3-5). \qed

\enddemo

\heading{\bf \S 4. Absolute torsion and Conway polynomial}\endheading

In this section we prove
that the absolute torsion of a 3-dimensio\-nal manifold $X=X_K$,
obtained by performing 
0-surgery on a knot $K\subset S^3$, is precisely {\it the Conway polynomial} of $K$. 
This suggests to view the absolute torsion as a generalization of the
Conway polynomial, which is applicable to high dimensions as well.

\subheading{4.1} Recall that the Conway link  polynomial is a function $L\mapsto \nabla_L$
from the set of isotopy classes of oriented links in $S^3$ into the ring
 of polynomials
$\Z [z]$. This function is uniquely
characterized by the following two properties:

(i) its value on a trivial knot is equal to 1;

(ii) for any skein triple of oriented links $L_+,L_-,L_0$,
$$\nabla_{L_+}(z)-\nabla_{L_-}(z)=z\nabla_{L_0}(z).\tag4-1$$
Here by a skein triple of links we mean three oriented links coinciding
outside a 3-ball and looking as the standard triple
(positive crossing of 2 strands, negative crossing of 2 strands, two
 vertical
strands) inside this ball.

Recall that the Conway polynomial $\nabla_K(z)$  of any oriented knot
$K$   involves only even powers
of $z$.

\subheading{4.2} Let $K\subset S^3$ be an oriented knot in $S^3$ and let the 3-manifold 
$X=X_K$ be obtained
by 0-surgery on $K$. Then $H_1(X)=\Z$ (has no 2-torsion) and the semi-characteristic $s\chi(X)=2$ is even. Also, the condition 3.1 is satisfied (since any orientable 3-manifold has a trivial 
tangent bundle). Hence, the absolute torsion
$\T(F)$ is well defined for arbitrary flat bundle $F$ over $X$ (cf. \S 3). It is independent of the
orientation of $X$. 

We will consider {\it line} flat bundles $F$ over $X$. 
It is clear that each such bundle is completely determined by the monodromy 
$a\in \C^\ast$ along a fixed generator $m\in H_1(X)$ (the meridian). We will denote this flat
line bundle by $F_a$.

It is well known that the homology $H_\ast(X;F_a)$ is trivial if
and only if $a\ne 1$ and $a$ is not
a root of the Alexander polynomial
of $K$. This excludes finitely many points of $\C^\ast$. On the complement of these points the
absolute torsion $\T(F_a)$
is a well defined $\C^\ast$-valued function of $a$.
We shall compute this function in terms of the  Conway polynomial of $K$.

\proclaim{4.3. Theorem} Let  $K\subset S^3$ be an oriented knot and $a\in \C^\ast$, $a\ne 1$,
is not a root of the Alexander polynomial of $K$. Then the absolute
torsion $\T(F_a)$ of the flat line bundle $F_a$ over $X=X_K$ 
is given by 
$$\T(F_a) \, =\, {\nabla_K(a^{1/2}
-a^{-1/2})}.\tag4-2$$ \endproclaim
\demo{Proof}
We begin by recalling the definition of $\nabla_K(z)$ in terms
of torsions, given in \cite{T1},
section 4.3. Let $Y$ be the exterior of $K$, i.e., the complement of
an open regular neighborhood of $K$. We provide $Y$ with a homology
orientation defined by the basis $([pt], t)$ where $[pt]$ is the
homology class of a point and $t$ is the   generator of
$H_1(Y)=\Z$ represented by a meridian of $K$. (Note that $H_i(Y)=0$ for
$i\neq 0,1$). Consider the (refined) Reidemeister torsion
  $\tau_0(Y)$ corresponding to the natural embedding of the ring
$\Z[H_1(Y)] =\Z[t,t^{-1}]$
into its field of fractions
$\Q(t)$. This torsion is an element of $\Q(t)$
 defined up to multiplication by
powers of $t$. Choose a representative $A(t)\in  \Q(t)$
of $\tau_0(Y)$. It is known that
$\overline {A(t)}=-t^m A(t)$ where $m\in \Z$ and the bar denotes
the involution in $\Q(t)$ sending $t$ to $t^{-1}$. Then
$$\nabla_K(t-t^{-1})=-(t-t^{-1})^{-1} t^m A(t^2).\tag4-3$$
Note that here we use a normalization of $\nabla$ different from the one
in \cite{T1}; this difference is responsible for the factor
$(t-t^{-1})^{-1}$ on the right hand side. Note also that the product
$ (1-t) A(t)$ is a polynomial in $t$ representing the Alexander
polynomial of $K$.

We can use the multiplicativity of torsions to compute the (refined)
Reidemeister torsion $\tau_0(X)$ corresponding to the natural embedding
of the ring  $\Z[H_1(X)]=\Z[H_1(Y)]
=\Z[t,t^{-1}]$
into its field of fractions
$\Q(t)$. Since we are dealing with the sign-refined torsions we need
to use the corresponding sign-refined multiplicativity theorem,
\cite{T1},
Theorem 3.4.1. By this theorem,
$\tau_0(X)=(-1)^\mu\tau_0(Y)\tau_0(X,Y)$
where $\mu$ is a certain residue modulo 2 and the pair $(X,Y)$ is
provided with a homology orientation induced by those in $X$ and $Y$,
as described in \cite{T1},
Theorem 3.4.1.  A direct computation shows that in our setting $\mu=0$
and $\tau_0(X,Y)=(1-t)^{-1}$. Thus,
$\tau_0(X)=(1-t)^{-1} \tau_0(Y)$.

For any Euler structure $\xi$ on $X$, we have a refined version
$\tau_0(X,\xi)\in \Q(t)$ of $\tau_0(X)$. By duality,
$\overline {\tau_0(X,\xi)}=c(\xi) \tau_0(X,\xi)$,
see \cite{T3},
section 2.7.
We take $\xi$ to be the canonical Euler structure on $X$,
so that
$\overline {\tau_0(X,\xi)}=  \tau_0(X,\xi)$.
By the argument above, the product
$A(t)=(1-t)\tau_0(X,\xi)$ is a representative of
$\tau_0(Y)$. It is clear that
$\overline {A(t)}=-t^{-1} A(t)$ so that
$$\nabla_K(t-t^{-1})=-(t-t^{-1})^{-1} t^{-1} A(t^2)=\tau_0(X,\xi)
(t^2)$$
where $\tau_0(X,\xi)
(t^2)$ is obtained from the rational function
$\tau_0(X,\xi)=\tau_0(X,\xi)(t)\in \Q(t)$ by doubling all powers of
$t$. It follows from definitions that
for any non-zero complex number  $a$, the number
$\T(F_a)$ is obtained from the rational function
$\tau_0(X,\xi)\in \Q(t)$ by the substitution $t=a$.
This implies the claim of the theorem. \qed
\enddemo

\subheading{4.4. Remark} There exist analogues of Theorem 4.3 for links, 
which involve the one-variable and multi-variable Conway polynomials. 

Let us briefly describe the case of one-variable Conway polynomial. 

Let $L=\{\ell_1, \ell_2,\dots, \ell_\mu\subset S^3\}$ be an oriented link. We will assume that 
the Conway polynomial of $L$ is nonzero.
Let $\lambda_{i, j}$ denote
the linking number of $i$-th and $j$-th components, $i\ne j$. Set 
$\lambda_j \, =\,  \sum \lambda_{i,j}$, where the summation is taken with respect to $i=1, \dots, \mu$,
$i\ne j$. Consider the closed 3-manifold $X=X_L$ obtained from $S^3$ by a surgery along the
link $L$, with framing $-\lambda_j$ along the component $\ell_j$, for
each $j=1, 2, \dots, \mu$. There is a canonical homomorphism $H_1(X)\to \Z$ (determined by the 
Seifert surface of $L$) and so, as above, for any complex number $a\in \C^\ast$ we have a canonical
flat line bundle $F_a$ over $X$. If $a\ne 1$ is not a root of the Alexander polynomial, then $H_\ast(X;F_a)=0$, and we obtain
$\T(F_a)\in \C^\ast$.
In this situation the following formula holds
$$\T(F_a) \, =\, \pm \, \, \frac{\nabla_L(a^{1/2}-a^{-1/2})}{(a^{1/2}-a^{-1/2})^{\mu-1}}.\tag4-4$$
It expresses the absolute torsion $\T(F_a)$
in terms of the Conway polynomial $\nabla_L(z)$ of $L$. 
Formula (4-4) allows also to find the 
Conway polynomial $\nabla_L(z)$ knowing $\T(F_a)$. To see this we recall that the Conway polynomial $\nabla_L(z)$ involves only powers $z^k$ with $k$ being of the same parity as 
$\mu-1$; hence $\nabla_L(z)/z^{\mu-1}$ is a Laurent polynomial in $z^2$.

\enddocument

\enddocument